\documentclass[twoside]{article}
\usepackage{amsfonts}
\usepackage{amsmath,amssymb}
\usepackage{graphicx}
\usepackage[mathlines]{lineno}

\newcommand*\patchAmsMathEnvironmentForLineno[1]{  \expandafter \let \csname old#1\expandafter \endcsname \csname
#1\endcsname
  \expandafter \let \csname oldend#1\expandafter \endcsname \csname
end#1\endcsname
  \renewenvironment{#1}     {\linenomath \csname old#1\endcsname}     {\csname oldend#1\endcsname \endlinenomath}}
\newcommand*\patchBothAmsMathEnvironmentsForLineno[1]{  \patchAmsMathEnvironmentForLineno{#1}  \patchAmsMathEnvironmentForLineno{#1*}}
\patchBothAmsMathEnvironmentsForLineno{equation}
\patchBothAmsMathEnvironmentsForLineno{align}
\patchBothAmsMathEnvironmentsForLineno{flalign}
\patchBothAmsMathEnvironmentsForLineno{alignat}
\patchBothAmsMathEnvironmentsForLineno{gather}
\patchBothAmsMathEnvironmentsForLineno{multline}

\numberwithin{equation}{section}

\begin{document}

\title{\vspace{-1in}%
\parbox{\linewidth}{\footnotesize\noindent
\textbf{Applied Mathematics E-Notes, ??(20??), ???-???} \copyright \hfill ISSN 1607-2510
\newline
Available free at mirror sites of http://www.math.nthu.edu.tw/$\sim$amen/} 
\vspace{0pt} \\
Sharp Bounds for the Arc Lemniscate Sine  Function\thanks{%
Mathematics Subject Classifications: 11M35, 26D07, 33B15, 33E20.}}
\date{{\small Received ????}}

\author{Horst Alzer\thanks{%
Morsbacher Stra{\ss}e 10, 51545 Waldbr\"ol, Germany,
\emph{email:} \tt{h.alzer@gmx.de}}\ , Man Kam Kwong\thanks{%
Department of Applied Mathematics, The Hong Kong Polytechnic University,
Hunghom, Hong Kong,
\emph{email:} \tt{mankwong@connect.polyu.hk
}}}
\maketitle

\begin{abstract}
The  arc lemniscate sine function is given by
$$
\mbox{arcsl}(x)=\int_0^x \frac{1}{\sqrt{1-t^4}}dt.
$$
In 2017, Mahmoud and Agarwal presented bounds for $\mbox{arcsl}$ in terms of the Lerch zeta function
$$
\Phi(z,s,a)=\sum_{k=0}^\infty \frac {z^k}{(k+a)^s}.
$$
They proved
$$
\frac{1}{8} \,  x \, \Phi(x^4, 3/2, 1/4) < \mbox{arcsl}(x)< \frac{1}{4} \, x  \, \Phi(x^4,3/2,1/4)\qquad{(0<x<1)}.
$$
We
 show that the factor $1/4$ can be replaced by $\mbox{arcsl}(1)/\Phi(1,3/2,1/4)=0.12836...$. This constant is best possible. 
\end{abstract}

\vspace{0.1cm}

\section{Introduction and statement of  result}

Let $F_1$ and $F_2$ be two points in the plane, with distance $F_1F_2=2c$. The lemniscate of Bernoulli is the locus of all points $ P $ such that  $PF_1 \cdot PF_2=c^2$. It is named after the Swiss mathematician Jakob Bernoulli (1655-1705) who was the first who studied the lemniscate in detail.
The arc length of the lemniscate curve $L$ is given by the formula
$$
L=4\sqrt{2} c \,  \mbox{arcsl}(1),
$$
where $\mbox{arcsl}$ is the so-called arc lemniscate sine function, defined by
$$
\mbox{arcsl}(x)=\int_0^x \frac{1}{\sqrt{1-t^4}} dt \qquad{(-1\leq x\leq 1)}.
$$
Many interesting information on this subject including  historical comments can be found in Ayoub \cite{A} and Langer \& Singer \cite{LS}.

This note is inspired by a remarkable paper published by
Mahmoud and Agarwal [1]  in 2017. Among others, the authors offered upper and lower bounds for $\mbox{arcsl}$
in terms of the
 Lerch zeta function
$$
\Phi(z,s,a)=\sum_{k=0}^\infty \frac{z^k}{(k+a)^s}.
$$
They proved the elegant double-inequality
\begin{equation}
\frac{1}{8} \,  x \, \Phi(x^4, 3/2, 1/4) < \mbox{arcsl}(x)< \frac{1}{4} \, x  \, \Phi(x^4,3/2,1/4) \qquad  (0<x<1).  \label{11}
\end{equation}
It is natural to ask whether the constant factors $1/8$ and $1/4$ are sharp.
In this note, we refine the upper bound given in (\ref{11}). Indeed, 
 the constant $1/4$ can be replaced by a smaller number as the following theorem reveals.

\vspace{0.3cm}
{THEOREM.} {For all $x\in (0,1)$ we have}
\begin{equation}
\alpha \, x  \, \Phi(x^4, 3/2, 1/4) < \mbox{arcsl}(x)< \beta \, x \, \Phi(x^4,3/2,1/4)  \label{12}
\end{equation}
{with the best possible constant factors}
\begin{equation}
\alpha=\frac{1}{8} \quad\mbox{and} \quad \beta=\frac{\mbox{arcsl}(1)}{\Phi(1,3/2,1/4)}=0.12836... .  \label{13}
\end{equation}

\vspace{0.3cm}
In particular, we obtain that for all $x\in (0,1)$ the ratio $\mbox{arcsl}(x)/ \bigl(x \Phi(x^4,3/2,1/4)\bigr)$ lies between $1/8$ and $1/7$.
The constant $\beta$ can be expressed in terms of
 the Euler  beta function and the Hurwitz zeta function, respectively, which are given by
$$
B(x,y)=\int_0^1 t^{x-1} (1-t)^{y-1} dt \quad\mbox{and}
\quad
\zeta(s,a)=\Phi(1,s,a)=\sum_{k=0}^\infty \frac{1}{(k+a)^s}.
$$
The substitution $t=s^{1/4}$ gives
$$
\mbox{arcsl}(1)=\int_0^1 \frac{1}{\sqrt{1-t^4}}dt=\frac{1}{4}\int_0^1 s^{-3/4} (1-s)^{-1/2}ds=\frac{1}{4} B(1/4,1/2).
$$
Thus,
$$
\beta=\frac{B(1/4,1/2)}{4 \zeta(3/2, 1/4)} = \frac{\Gamma(1/4)^2}{4\sqrt{2\pi} \zeta(3/2, 1/4)} \, ,
$$
where $\Gamma$ denotes the classical gamma function.

Schneider \cite{S} proved in 1937 that the lemniscate constant $\mbox{arcsl}(1)$ is a transcendental number; see also Todd~\cite{T}.

\vspace{0.4cm}
\section{Proof of  Theorem}

The following lemma plays an important role in the proof of our theorem. It is known in the literature as the monotone form of l'Hopital's rule; see Hardy et al. \cite[p.\ 106]{HLP}, Kwong \cite{K} and Pinelis \cite{P}.

\vspace{0.3cm}
{LEMMA.} {Let $u,v: [a,b]\rightarrow\mathbb{R}$ be continuous functions. Moreover, let $u,v$ be differentiable on $(a,b)$ and $v'\neq 0$ on $(a,b)$. If $u'/v'$ is strictly increasing on $(a,b)$, then
$$
x\mapsto \frac{u(x)-u(a)}{v(x)-v(a)}
$$
is strictly increasing on $(a,b)$.}

\vspace{0.4cm}
PROOF OF THEOREM. Let
$$
F(x)=\frac{\mbox{arcsl}(x)}{x \Phi(x^4, 3/2, 1/4)}.
$$
In order to prove that $F$ is strictly increasing on $(0,1)$ we
 apply the lemma with
$$
u(x)=\mbox{arcsl}(x) \quad\mbox{and} \quad v(x)=x\Phi(x^4, 3/2, 1/4).
$$
Let $x\in (0,1)$. We have
$$
u(0)=v(0)=0\quad\mbox{and} \quad u'(x)=\frac{1}{\sqrt{1-x^4}}, \quad v'(x)=8\sum_{k=0}^\infty \frac{x^{4k}}{\sqrt{4k+1}}.
$$
It follows that
\begin{equation}
\frac{u'(x)}{v'(x)}=\frac{1}{8 h(x^4)}  \label{21}
\end{equation}
with
$$
h(s)=\sqrt{1-s}\sum_{k=0}^\infty \frac{s^k}{\sqrt{4k+1}}.
$$
Then,
$$
2\sqrt{1-s} h'(s)=2(1-s)\sum_{k=1}^\infty \frac{k s^{k-1}}{\sqrt{4k+1}}-\sum_{k=0}^\infty \frac{s^k}{\sqrt{4k+1}}
=\sum_{k=0}^\infty a_k s^k,
$$
where
$$
a_k=\frac{2k+2}{\sqrt{4k+5}}-\frac{2k+1}{\sqrt{4k+1}}
=\frac{-1}{(2k+2)(4k+1)\sqrt{4k+5} +(2k+1)(4k+5)\sqrt{4k+1}   }.
$$
Since $a_k<0$ for $k=0,1,2,...$, we conclude that $h'(s)<0$ for $s\in (0,1)$. Thus, $h$ is strictly decreasing on $(0,1)$. Using
(\ref{21}) yields that $u'/v'$ is strictly increasing on $(0,1)$, so that the lemma reveals that $F=u/v$ is strictly increasing on $(0,1)$. It follows that
\begin{equation}
F(0) < F(x) < F(1)  \quad\mbox{for} \quad x\in (0,1).  \label{22}
\end{equation}
We have
\begin{equation}
F(1)=
\frac{\mbox{arcsl}(1)}{\Phi(1, 3/2, 1/4)}.  \label{23}
\end{equation}
Since
$$
\lim_{x\to 0} \frac{\mbox{arcsl}(x)}{x}
=\left . \frac{d}{dx}\mbox{arcsl}(x) \, \right  |_{x=0}=1
\quad\mbox{and} \quad \Phi(0,3/2,1/4)=8,
$$
we obtain
\begin{equation}
F(0)=\frac{1}{8}.  \label{24}
\end{equation}
From (\ref{22}), (\ref{23}) and (\ref{24}) we conclude that (\ref{12}) is valid and that the
  constant factors  $\alpha$ and $\beta$ as given in (\ref{13}) are best possible.
This completes the proof of the theorem.
\def\BOX #1 #2 {\framebox[#1in]{\parbox{#1in}{\vspace{#2in}}}}
\hfill \raisebox{1pt}{\BOX 0.1  -0.03 }


\vspace{1cm}
\begin{thebibliography}{99}

\bibitem{A}
R. Ayoub, The lemniscate and Fagnano's contributions to elliptic integrals, Arch. Hist. Exact Sci. 29 (1984), 131-149.

\bibitem{HLP}
G.H. Hardy, J.E. Littlewood, G. P\' olya, Inequalities, Camb. Univ. Press, Cambridge, 1952.


\bibitem{K}
M.K. Kwong, On Hopital-style rules for monotonicity and oscillation, arXiv:1502.07805 [math.CA] (2015).

\bibitem{LS}
J.C. Langer, D.A. Singer, Reflections on the lemniscate of Bernoulli: The forty-eight faces of a mathematical gem, Milan J. Math. 78 (2010), 643-682.

\bibitem{MA}
M. Mahmoud, R.P. Agarwal, On some bounds of Gauss arc lemniscate sine and tangent functions, J. Inequal. Spec. Func. 8 (2017), 46-58.


\bibitem{P}
I. Pinelis, L'Hospital type rules for monotonicity, with applications, J. Inequal. Pure Appl. Math. 3 (2002), article 5, 5 pp.

\bibitem{S}
T. Schneider, Arithmetische Untersuchungen elliptischer Integrale, Math. Ann. 113 (1937), 1-13.

\bibitem{T}
J. Todd, The lemniscate constants, Comm. ACM 18 (1975), 14-18.

\end{thebibliography}
\end{document}